\documentclass[preprint,12pt]{elsarticle}
\usepackage{verbatim}
\usepackage{amssymb}
\usepackage{amsmath,amssymb,amsfonts,graphicx,float,multicol,fleqn}
\usepackage{multirow}

\journal{.}

\begin{document}
\begin{frontmatter}
\title{The meshless method for solving radiative transfer problems in a slab medium based on radial basis functions}

\author[a]{J. A. Rad\corref{cor}}
\ead{j.amanirad@gmail.com}
\author[b]{S. Kazem}
\ead{saeedkazem@gmail.com}
\author[a]{K. Parand}
\cortext[cor]{Corresponding author. Tel:+98 937 509 2738.}
\ead{k\_parand@sbu.ac.ir}

\address[a]{Department of Computer Sciences, Shahid Beheshti University, G.C., Tehran, Iran}
\address[b]{Department of Mathematics, Imam Khomeini International University, Ghazvin 34149-16818, Iran}

\begin{abstract}
In this paper a numerical meshless method for solving the radiative transfer equations in a slab medium with an isotropic scattering is considered. The method is based on radial basis functions to approximate the solution of an integral-partial differential equation by using collocation method.
For this purpose different applications of RBFs are used. To this end the numerical solutions are obtained without any mesh generation into the domain of the problems.
The results of numerical experiments are compared with the existing results in illustrative examples to confirm the
accuracy and efficiency of the presented scheme. Also the norm of the residual functions are obtained to show the convergence of the method.
\end{abstract}
        \begin{keyword}
        Radiative transfer equation, Radial basis functions, Integral-partial differential equation.
        \end{keyword}
    \end{frontmatter}

\section{Introduction}
\subsection{Radiative transfer equations}
The radiative transfer equations (RTE) contain comprehensive applications in lightweight fibrous insulation, the study of atmospheres, coal-fired combustion and conversion systems, and in remote sensing \cite{radiative1,radiative3}. Solving these types of equations have been of considerable interest.
In this paper, the RTE in slab medium is considered. It is irradiated by an isotropic radiation field $I_0$, with azimuthal symmetry, finally this equation has been obtained in references \cite{radiative1,radiative2} as:
\begin{eqnarray}\label{radiative}
\frac{x}{t_0}\frac{\partial}{\partial y}I(y,x)+I(y,x)=S(y)+\frac{\omega}{2}\int_{-1}^{1}P(x,\hat{x})I(y,\hat{x})d\hat{x}~,
\end{eqnarray}
where $I$ is the angle distribution of intensity normalized to $I_0$, $y$ the distance normalized by the optical depth $t_0$ of the slab, $x$ the direction cosine of the angle made by the specific intensity at any depth $y$ with the direction of increasing $y$, $\omega$ the albedo of a single scattering and $S(y)$ the dimensionless emission source. The boundary conditions of the problem are:
\begin{eqnarray}\label{b1}
&&I(0,x)=I^0(x)~,~~~~~~0<x\leq1~,\\\label{b2}
&&I(1,x)=I^1(x)~,~~~~-1\leq x<0~,
\end{eqnarray}
functions $I^0(x)$ and $I^1(x)$ are known. The phase scattering function $P(x,\hat{x})$, is represented in term
of the Legendre polynomials of the first kind $P_n(x)$ by the expansion \cite{radiative1,radiative3}
\begin{eqnarray}\label{PP}
&&P(x,\hat{x})=\sum_{i=0}^{n}c_iP_i(x)P_i(\hat{x})~,\\\nonumber
&&P_0(x)=1~,\qquad P_1(x)=x~,\\\nonumber
&&P_{n+1}(x)=\frac{2n+1}{n+1}xP_n(x)-\frac{n}{n+1}P_{n-1}(x)~,~~~n\geq1
\end{eqnarray}
where $c_i$ are the expansion coefficients with $c_0=1$. The RTE in slab medium with anisotropic scattering has some numerical and rigorous solutions such as: two-flux \cite{radiative6,radiative7}, spherical harmonic \cite{radiative5,radiative8}, series expansion \cite{radiative9,radiative10,radiative11}, integral equation \cite{radiative14,radiative15}, Pad\'{e} approximation \cite{radiative18}, iterative \cite{radiative5}, variational \cite{radiative14,radiative12,radiative13}, eigenfunction expansion \cite{radiative16,radiative17}, the linear spline approximation \cite{radiative1}, the generalized Eddington approximation \cite{radiative2} and Spectral methods approximation \cite{TAU,galerkin,collocation}.\\
In this paper, a new approach to the solution of RTE is presented. Our approach is based on radial basis functions (RBFs) collocation method to approximate unknown function $I(y,x)$ and solve the Eq. (\ref{radiative}) with Eqs. (\ref{b1}), (\ref{b2}) and (\ref{PP}).

\subsection{Radial basis functions}
RBFs interpolation are techniques for representing a function starting with data on scattered nodes. This technique first appears in the literature as a method for scattered data interpolation, and the method was highly favored after being reviewed by Franke \cite{Franke}, who found it to be the most impressive of the many methods he tested. Later, Kansa \cite{Kansa.127.145,Kansa.EJ1990} proposed a scheme for the estimation of
partial derivatives using RBFs. The main advantage of radial basis functions methods is the
meshless characteristic of them. The use of radial basis functions as a meshless method for the numerical solution
of partial differential equations (PDEs) is based on the collocation method. These methods have recently received a great deal of attention from researchers \cite{M.sharan,Zerroukat.M1998,Mai-Duy.N2001,tatari.dehghan.2010,Dehghan.A.Shokri,A.Alipanah.M.Dehghan,Sarra}.

Recently, RBFs methods were extended to solve various ordinary and partial differential equations including the high order ordinary differential equations \cite{N.Mai-Duy2005}, second-order parabolic equation with nonlocal boundary conditions \cite{tatari.dehghan2007}, the nonlinear Klein-Gordon equation \cite{Dehghan.A.Shokri}, regularized long wave (RLW) equation \cite{ali.haqb.Islam}, Hirota-Satsuma coupled KdV equations \cite{Khattak.Tirmizi.Islam}, a system of nonlinear integral equations \cite{Golbabai2009}, Second-order hyperbolic telegraph equation \cite{MehdiDehghan.ArezouGhesmati}, the solution of 2D biharmonic equations \cite{Mai-Duy2009}, the case of heat transfer equations \cite{kazem.darcian} and so on.
\\
One of the most powerful interpolation method with analytic two-dimensional test function is the RBFs method based on multiquadric (MQ) basis function
\begin{align}\label{MQ}
\phi(r)=\sqrt{r^2+c^2}~,
\end{align}
suggested by R.L. Hardy \cite{Hardy}. Madych and Nelson \cite{madych} showed that interpolation with MQ is exponentially convergent based on reproducing kernel Hilbert space. Wu and Schaback \cite{Wu.1993} use a different technique to handle the case of interpolation with power spline and the thin plate spline. Convergence property of the MQ has been also showed by Buhman \cite{M.D.Buhmann,buhmann.book}. Too large or too small shape parameter $c$ in Eq. (\ref{MQ})  make the MQ too flat and too peaked. Despite many studies done to find algorithms for selecting the optimum values of $c$ \cite{S.Rippa,Cheng.Golberg.Kansa.Zammito,Carlson.Foley,A.E.Tarwater,G.E.Fasshauer.J.G.Zhang}, the optimal choice of shape parameter is an open problem which is still under intensive investigation.\\
For more basic details about RBFs the interested readers can refer to the recent books and paper by Buhmann \cite{M.D.Buhmann,buhmann.book} and Wendland \cite{HOLGER.WENDLAND}, compactly and globally supported; and convergence rate of the radial basis functions.

Some of the infinitely smooth RBFs choices are listed in Table \ref{Tab.RBR.definition}. The RBFs can be of various types, for example: multiquadrics (MQ), inverse multiquadrics (IMQ), Gaussian forms (GA) form etc.
Regarding the inverse quadratic, inverse
multiquadric (IMQ) and Gaussian (GA), the coefficient matrix of RBFs interpolating is positive definite and, for
multiquadric (MQ), it has one positive eigenvalue and the
remaining ones are all negative \cite{M.J.D.Powell}.\\
This paper is arranged as follows: in Section \ref{RCF.Intro}, we describe the properties of radial basis functions. In Section \ref{solution} we implement the problem with the proposed method , report our numerical finding and demonstrate the accuracy of the proposed methods. The conclusions are discussed in the final Section.
%*****************************************************************************************************************
\section{Radial basis functions}\label{RCF.Intro}
\subsection{Definition of radial basis functions}
Let $\mathbb{R}^{+}=\{x\in \mathbb{R},x\geq 0\}$ be the non-negative half-line and let $\phi:\mathbb{R}^{+}\to \mathbb{R}$ be a continuous function with $\phi(0)\geq 0$. A radial basis function on $\mathbb{R}^{d}$ is a function of the form
\begin{eqnarray}\nonumber
    \phi(\|X-X_i\|)~,
\end{eqnarray}
where $X, ~X_i \in \mathbb{R}^{d}$ and $\|.\|$ denotes the Euclidean distance between $X$ and $X_i$s. If one chooses $N$ points $\{X_i\}_{i=1}^{N}$ in $\mathbb{R}^d$ then by custom
\begin{eqnarray}\nonumber
s(X)=\sum_{i=1}^{N}\lambda_i\phi(\|X-X_i\|);\quad \lambda_i \in \mathbb{R}~,
\end{eqnarray}
is called a radial basis function as well \cite{M.A.Golberg}.\\
The standard radial basis functions are categorized into two major classes \cite{Khattak.Tirmizi.Islam}:\\
\textbf{Class 1}. Infinitely smooth RBFs \cite{Khattak.Tirmizi.Islam,Dehghan·Shokri}:\\
These basis functions are infinitely differentiable and heavily depend on the shape parameter $c$ e.g. Hardy multiquadric (MQ), Gaussian(GA), inverse multiquadric (IMQ), and inverse quadric(IQ)(See Table \ref{Tab.RBR.definition}).

\textbf{Class 2}. Infinitely smooth (except at centers) RBFs \cite{Khattak.Tirmizi.Islam,Dehghan·Shokri}:\\
The basis functions of this category are not infinitely differentiable. These basis functions are shape parameter free and have comparatively less accuracy than the basis functions discussed in the Class 1. For example, thin plate spline, etc \cite{Khattak.Tirmizi.Islam}.\\
%**************************************************************
\subsection{RBFs interpolation}
The $d$-dimensional function $F(X)$, $F:\mathbb{R}^d\longrightarrow\mathbb{R}$
, to be interpolated or approximated can be
represented by an RBFs as:
\begin{eqnarray}\label{rbf intrerpolation}
F(X)\approx F_N(X)=\sum_{i=1}^{N}\lambda_{i}\phi_{i}(r)=\mathbf{\Phi}^T(r)\mathbf{\Lambda}~,
\end{eqnarray}
where
\begin{eqnarray}\nonumber
&&\phi_i(r)=\phi(\|X-X_{i}\|)~,\\\nonumber
&&\mathbf{\Phi}^T(r)=[\phi_{1}(r), \phi_{2}(r), ..., \phi_{N}(r)],\\\label{lambda}
&&\mathbf{\Lambda}=[\lambda_{1}, \lambda_{2}, ..., \lambda_{N}]^T,
\end{eqnarray}
$X$ is the input and $\{\lambda_{i}\}_{i=1}^N$ are the set of coefficients to be determined. By choosing $N$ collocation nodes $\{X_{i}\}_{i=1}^N$, we can approximate the function $F(X)$.
\begin{eqnarray}\nonumber
y_j=F(X_j)=\sum_{i=1}^{N}\lambda_{i}\phi_{i}(r_{j})  , \quad j=1, 2, ..., N~,
\end{eqnarray}
where $r_j=\|X_j-X_{i}\|$. To summarize discussion on coefficient matrix, we define:
\begin{eqnarray}\nonumber
\mathbf{A}\mathbf{\Lambda}=\mathbf{Y}~,
\end{eqnarray}
where
\begin{align}\nonumber
\mathbf{Y}=&[y_{1}, y_{2}, ..., y_{N}]^T,\\\nonumber
\mathbf{A}=&[\mathbf{\Phi}^T(r_{1}), \mathbf{\Phi}^T(r_{2}), ..., \mathbf{\Phi}^T(r_{N})]^T~,\\
 =&
\left[
 \begin{array}{cccc}\label{matrice}
\phi_{1}(r_{1}) & \phi_{2}(r_{1}) & \dots  & \phi_{N}(r_{1}) \\
\phi_{1}(r_{2}) & \phi_{2}(r_{2}) & \dots  & \phi_{N}(r_{2}) \\
\vdots          & \vdots          & \ddots & \vdots          \\
\phi_{1}(r_{N}) & \phi_{2}(r_{N}) & \dots  & \phi_{N}(r_{N}) \\
\end{array}
\right]~.
\end{align}
We have $\phi_{i}(r_{j})=\phi_{j}(r_{i})$ consequently $\mathbf{A}=\mathbf{A}^T$.\\
All the infinitely smooth RBFs choices are listed in Table \ref{Tab.RBR.definition} will give
coefficient matrices $\mathbf{A}$ in Eq. (\ref{matrice}) which are symmetric and nonsingular \cite{M.J.D.Powell},
i.e. there is a unique interpolant of the form Eq. (\ref{rbf intrerpolation})
no matter how the distinct data points are scattered in any number
of space dimensions. \\
%********************************************************************************************
\section{RBFs collocation method for solving RTE}\label{solution}
\begin{figure}
\centerline{\includegraphics[scale=0.5]{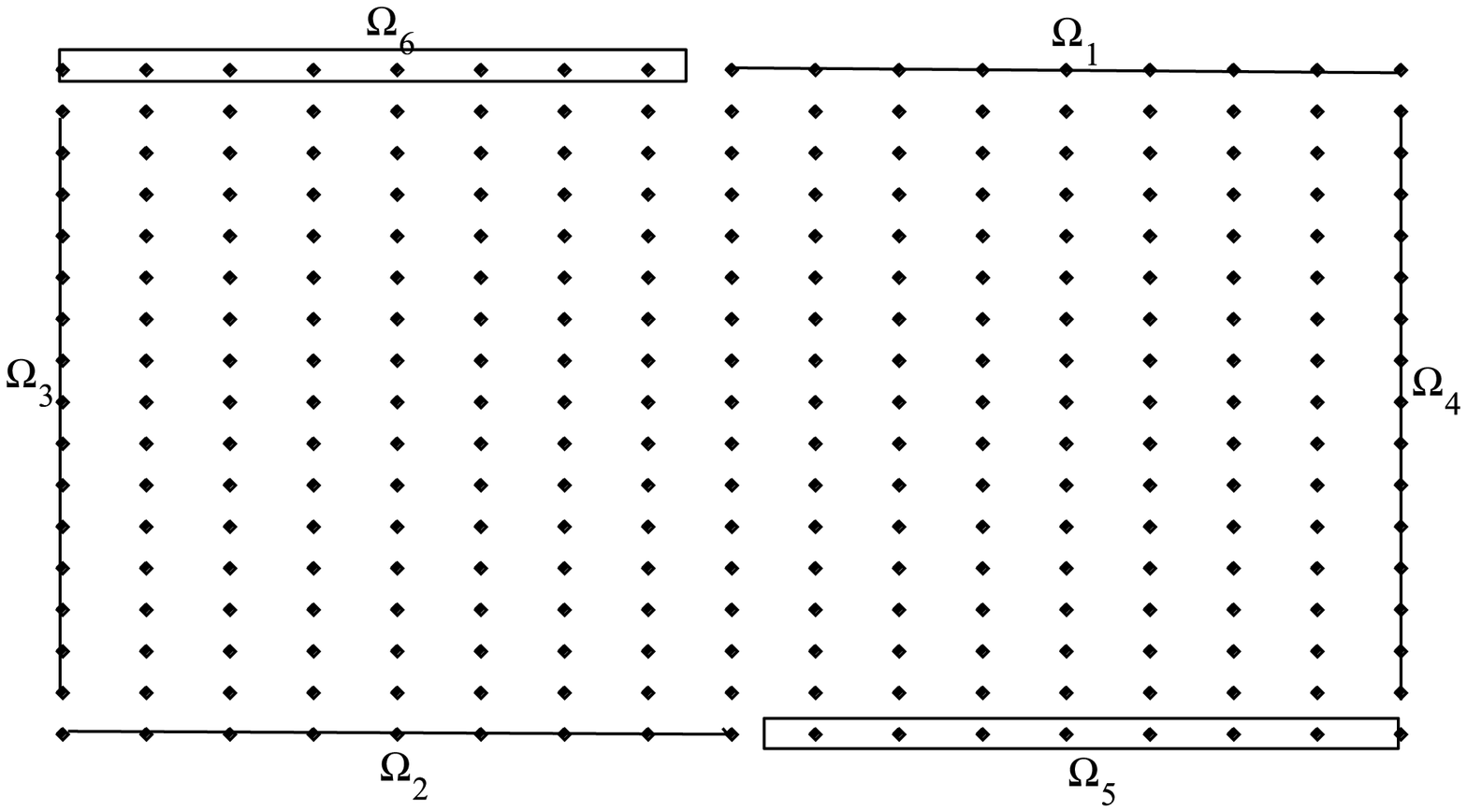}}
\caption{Graph of center nodes $\{X_k\}_{k=1}^{N}={(y_i,x_j)}$ }
\label{fig1}
\end{figure}

We approximate 2-dimensional $I(y,x)$ by RBFs as:
\begin{eqnarray}\label{I(y,x)}
&&I(y,x)\simeq I_N(y,x)=\sum_{k=1}^{N}\lambda_k\phi(\|X-X_{k}\|)~,\\\label{I^'(y,x)}
&&\frac{\partial}{\partial y}I(y,x)\simeq \frac{\partial}{\partial y}I_N(y,x)=\sum_{k=1}^{N}\lambda_k\frac{\partial}{\partial y}\phi(\|X-X_{k}\|)~,
\end{eqnarray}
where $\|X-X_{k}\|=\sqrt{(y-y_i)^2+(x-x_j)^2}$~, $-1\leq x\leq 1$~, $0\leq y\leq 1$ and $i=0,1,...,m$~, $j=0,1,...,n$ and $N=(n+1)(m+1)$.\\
Center nodes $(y_i,x_j)$ are chosen as:
\begin{eqnarray}
&&y_i=\frac{i}{m}~,~~~~~~~~\qquad i=0,1,...,m~, \\
&&x_j=\frac{(2j-n)}{2n}~,\qquad j=0,1,...,n~.
\end{eqnarray}
We define $Res(y,x)$ by substitute Eqs. (\ref{I(y,x)}), (\ref{I^'(y,x)}) in Eq. (\ref{radiative})
\begin{eqnarray}\nonumber
Res(y,x)=\frac{x}{t_0}~\frac{\partial}{\partial y}I_N(y,x)+I_N(y,x)-S(y)-\frac{\omega}{2}\int_{-1}^{1}P(x,\hat{x})I_N(y,\hat{x})d\hat{x}\\\label{Res}
\end{eqnarray}
The boundary conditions of the problem Eqs. (\ref{b1}), (\ref{b2}) are obtained as
\begin{eqnarray}\label{boundary}
\begin{cases}
I_N(0,x)\simeq I^0(x)~,~ 0<x \leq 1~,\\\cr
I_N(1,x)\simeq I^1(x)~,~-1\leq x<0~.
\end{cases}
\end{eqnarray}
Now to obtain the coefficients $\{\lambda_k\}_{k=1}^{N}$, we use $N$ collocation nodes the same as center nodes $\{X_k\}_{k=1}^{N}$ in Figure \ref{fig1}. A set of linear algebraic equations is constructed by discretizing Eq. (\ref{Res}) on five sets of equations, and using two boundary conditions appeared in Eq. (\ref{boundary}) as
\begin{eqnarray}\label{R}
\begin{cases}
Res(1,x_j)=0~,~\quad  ~~~~~~ \Omega_1=\{y=1~,~x_j~;~j=\frac{n}{2},\frac{n}{2}+1,...,n\}~,\\\cr
Res(0,x_j)=0~,\quad  ~~~~~~~ \Omega_2=\{y=0~,~x_j~;~j=0,1,...,\frac{n}{2}\}~,\\\cr
Res(y_i,-1)=0~,\quad ~~~~~ \Omega_3=\{x=-1~,~y_i~;~i=1,2,...,m-1\}~,\\\cr
Res(y_i,1)=0~,~~\quad ~~~~~ \Omega_4=\{x=1~,~y_i~;~i=1,2,...,m-1\}~,\\\cr
Res(y_i,x_j)=0~,~~\quad 0< y_i< 1~,~-1<x_j<1~;~i=1,2,...,m-1~,~n=1,2,...,n-1~,\\\cr
I_N(0,x_j)= I^0(x_j)~,~ ~~~~~\Omega_5=\{y=0~,~x_j~;~j=\frac{n}{2}+1,\frac{n}{2}+2,...,n\}~,\\\cr
I_N(1,x_j)= I^1(x_j)~,~~~~~~\Omega_6=\{y=1~,~x_j~;~j=0,1,...,\frac{n}{2}-1\}~.\\\cr
\end{cases}
\end{eqnarray}
The number of unknown coefficients $\{\lambda_k\}_{k=1}^{N}$ is equal to $(m+1)(n+1)$ and these can be
obtained from Eq. (\ref{R}). Consequently, $I(y,x)$ can be calculated.\\
For RTE the radiative fluxes
\begin{eqnarray}
F^+(y)=2\int_{0}^{1}I(y,x)xdx~,
\end{eqnarray}
at the lower boundary $F^+(1)$ is important. We can achieve $F^+(1)$ by using approximate function as
\begin{eqnarray}
F^+(1)\simeq2\int_{0}^{1}I_N(1,x)xdx~.
\end{eqnarray}
In some cases, exact value of $F^+(1)$ is reported in \cite{radiative23}.
Two cases of radiative transfer equations, are given in \cite{radiative1,radiative2,radiative23}, presented as test examples to show the reliability of the method.
%************************************************************************
\subsection{Example 1}
Let the RTE be given as \cite{radiative1,TAU,galerkin,collocation}
\begin{eqnarray}\label{exe1}
\begin{cases}
\frac{x}{t_0}\frac{\partial}{\partial y}I(y,x)+I(y,x)=\frac{1}{2}\int_{-1}^{1}[1+c_1P_1(x)P_1(\hat{x})]I(y,\hat{x})d\hat{x}~,\\\cr
I(0,x)=1~,~0<x\leq1~,~I(1,x)=0~,~-1\leq x<0~.
\end{cases}
\end{eqnarray}
We applied the present method and solved Eq. (\ref{exe1}) and then evaluated the radiative
fluxes at the lower boundary $F^+(1)$ for $c_1=0.7,0,-0.7$ and $t_0=0.1,~0.5,~1.0,~3.0$. In Table \ref{Tab2}, achievement values of $F^+(1)$ with $n=m=20$ are compared with exact values reported in \cite{radiative23}.
%*****************************************************************************
\subsection{Example 2}
This case of RTE was first considered by Menguc and Viskanta in \cite{radiative24}
\begin{eqnarray}\label{exe2}
\begin{cases}
x\frac{\partial}{\partial y}I(y,x)+I(y,x)=\frac{0.8}{2}\int_{-1}^{1}[1+\sum_{l=1}^{4}c_lP_l(x)P_l(\hat{x})]I(y,\hat{x})d\hat{x}~,\\\cr
I(0,x)=1~,~0<x\leq1~,~I(1,x)=0~,~-1\leq x<0~,
\end{cases}
\end{eqnarray}
where, $c_1=0.6438$, $c_2=0.5542$, $c_3=0.1036$, $c_4=0.0105$.\\
We applied the present method and solved Eq. (\ref{exe2}) and then evaluated $F^+(1)$. As a result, in Table \ref{Tab3}, value of $F^+(1)$ obtained by the present method with $n=m=20$, is compared with different approximate methods.
%***************************************************************************************************************************

%**********************************************************************************************************************
%*********************************************************************************************************************

%*******************************************Conclusions***********************************************************************************
\section{Concluding remark}\label{sec.conclusion}
\begin{figure}
\centering
\includegraphics[width=4in]{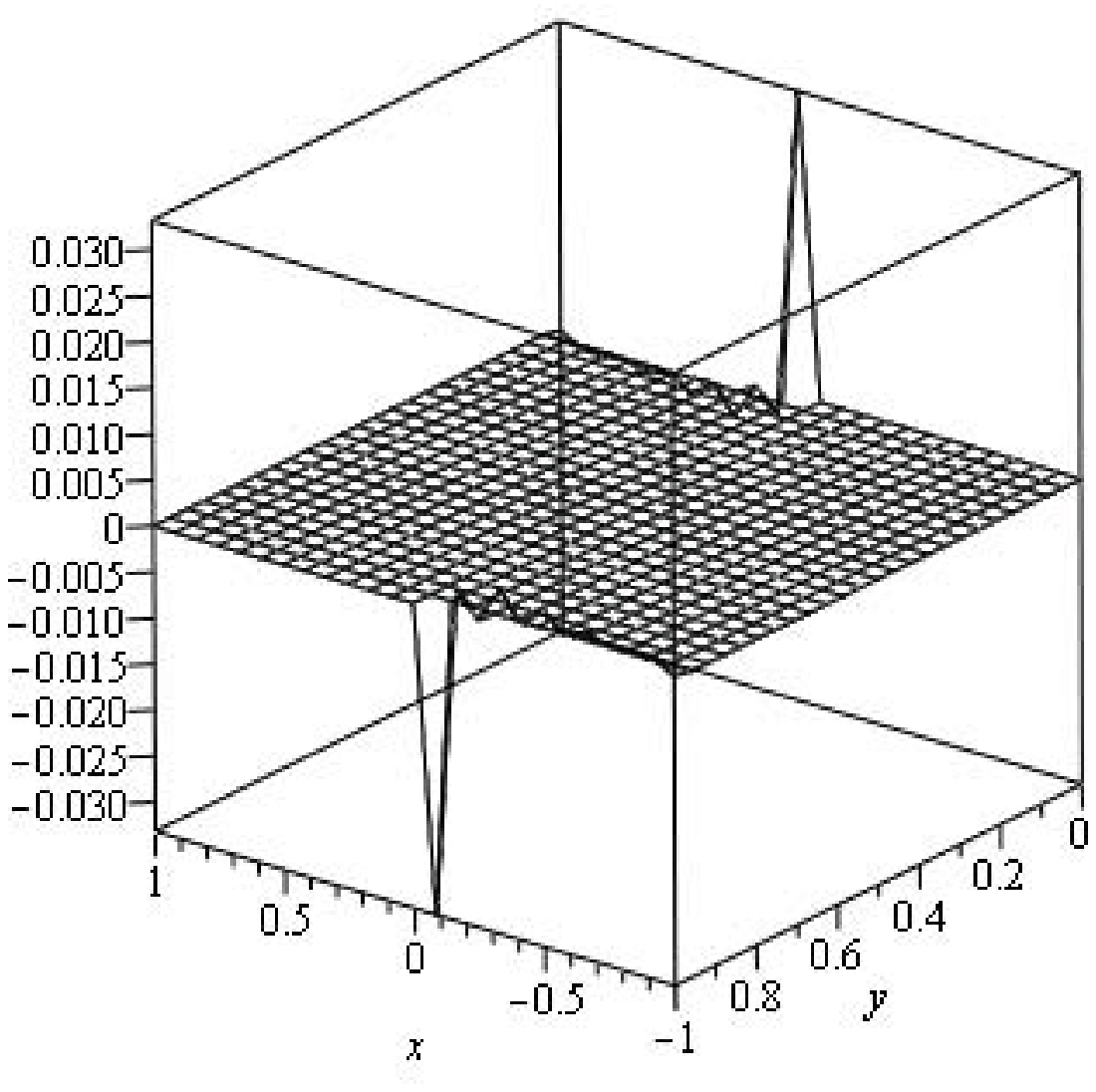}
\caption{Figure of $Res(y,x)$ Example 1 for case of $t_0=3$, $c_0=0.7$ by using MQ with $c=0.3$ and n=m=24.}
\label{Fig2}
\end{figure}

\begin{figure}
\centering
\includegraphics[width=4in]{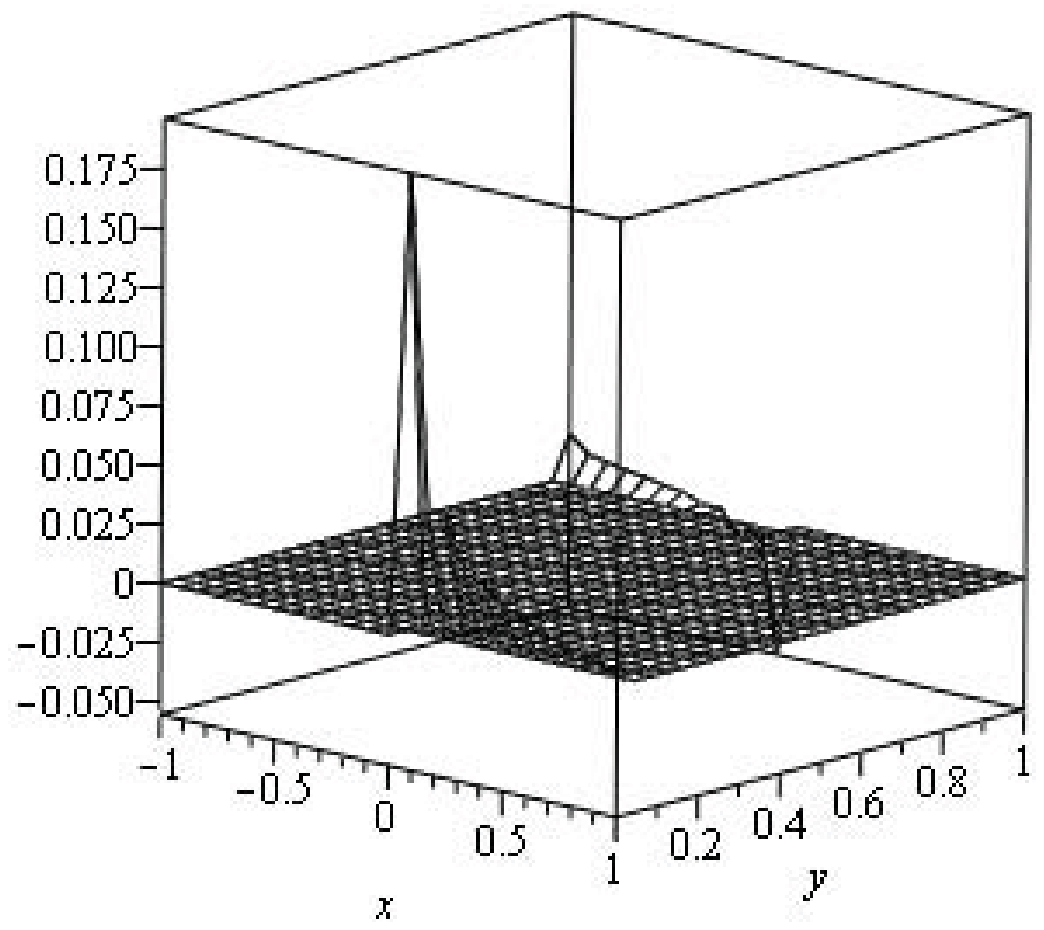}
\caption{Figure of $Res(y,x)$ Example 2 by using MQ with $c=0.3$ and n=m=24.}
\label{Fig3}
\end{figure}

\begin{figure}
\centering
\includegraphics[width=4in]{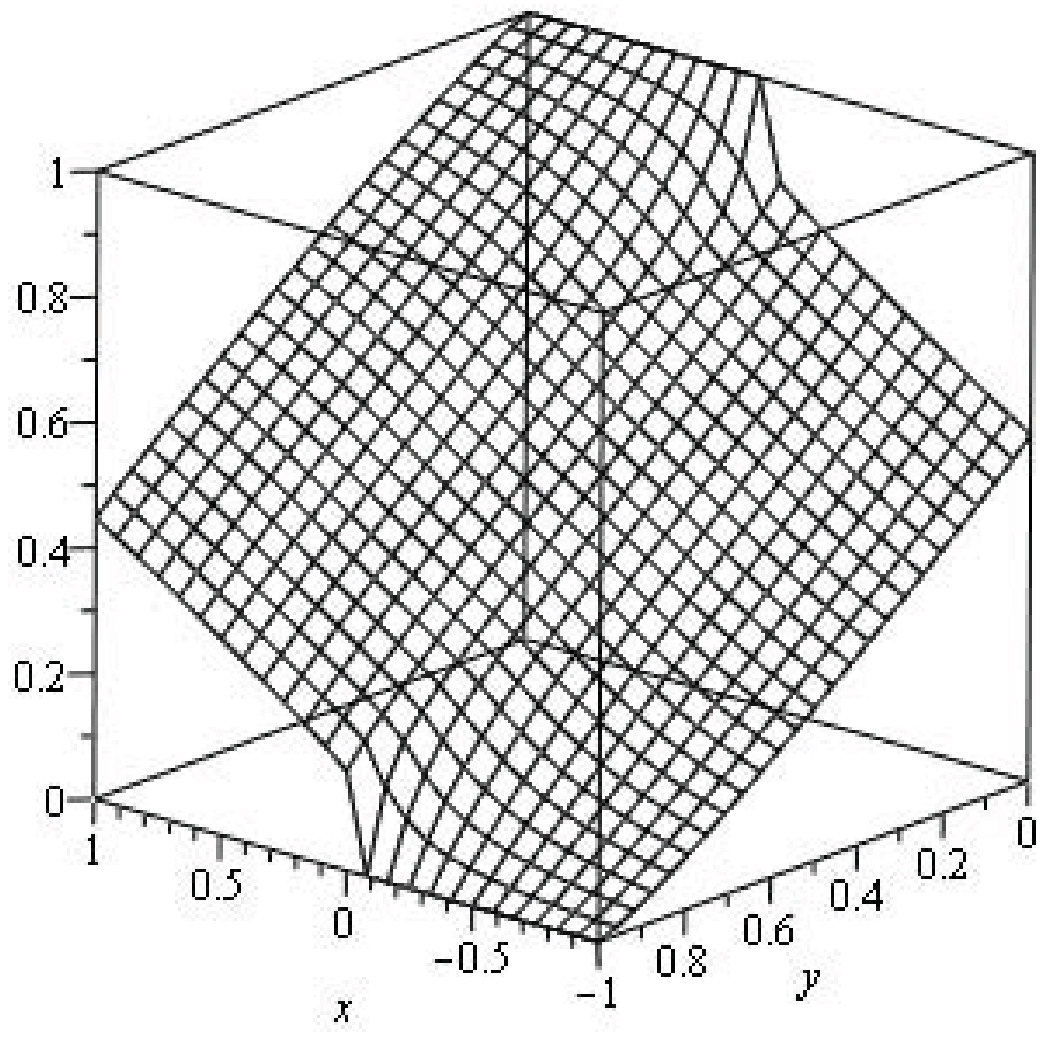}
\caption{Figure of $I(y,x)$ Example 1 for case of $t_0=3$, $c_0=0.7$ by using MQ with $c=0.3$ and n=m=24.}
\label{Fig4}
\end{figure}

\begin{figure}
\centering
\includegraphics[width=4in]{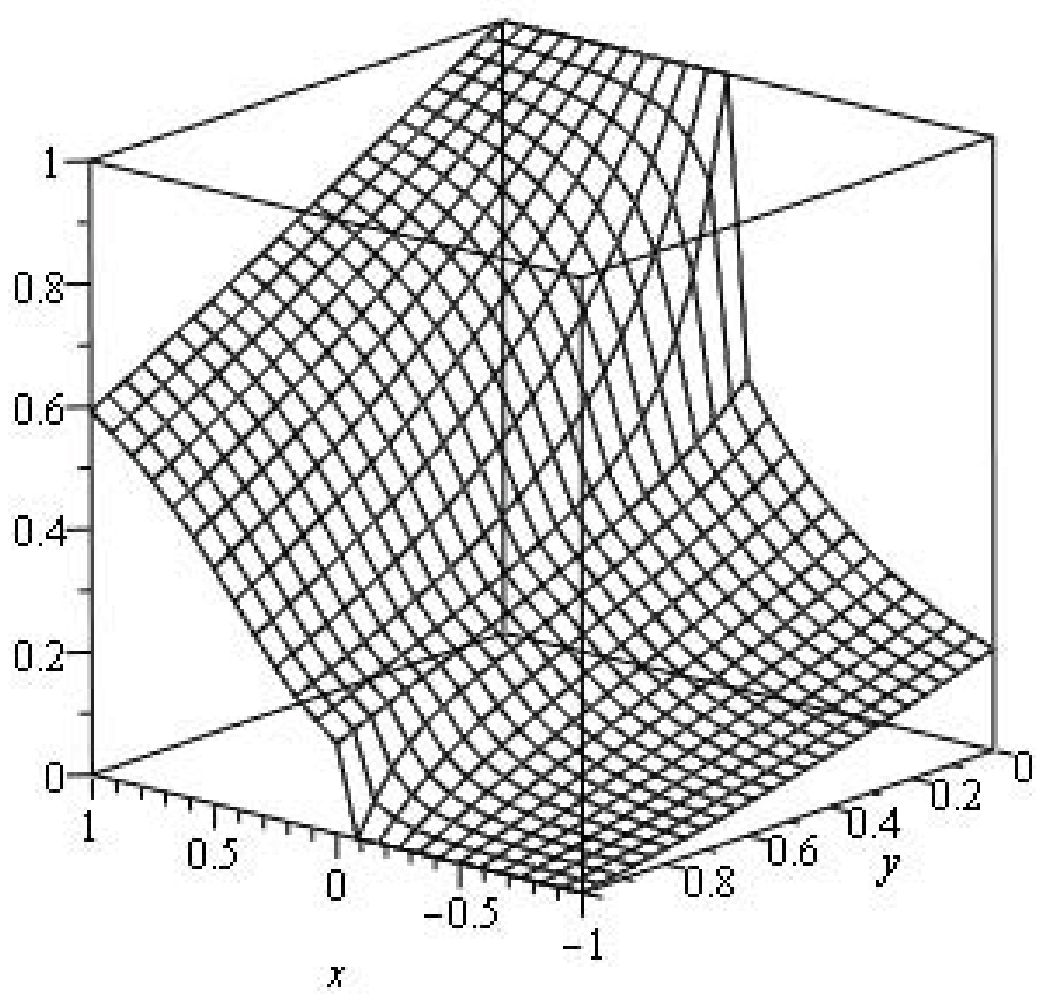}
\caption{Figure of $I(y,x)$ Example 2 by using MQ with $c=0.3$ and n=m=24.}
\label{Fig5}
\end{figure}

\begin{figure}
\centering
\includegraphics[width=3in]{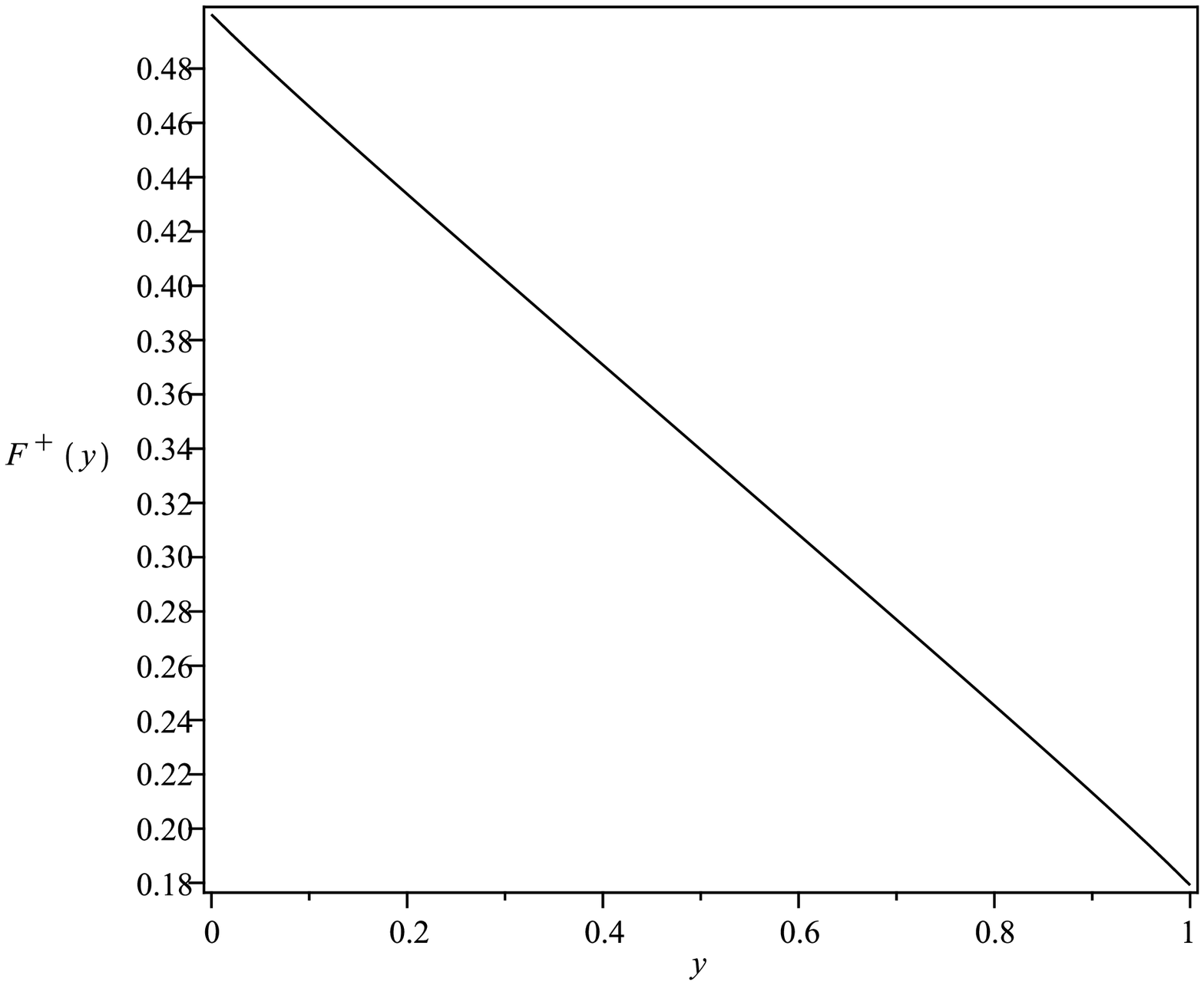}
\caption{Figure of $F^+(y)$ Example 1 for case of $t_0=3$, $c_0=0.7$ by using MQ with $c=0.3$ and n=m=24.}
\label{Fig6}
\end{figure}

\begin{figure}
\centering
\includegraphics[width=3in]{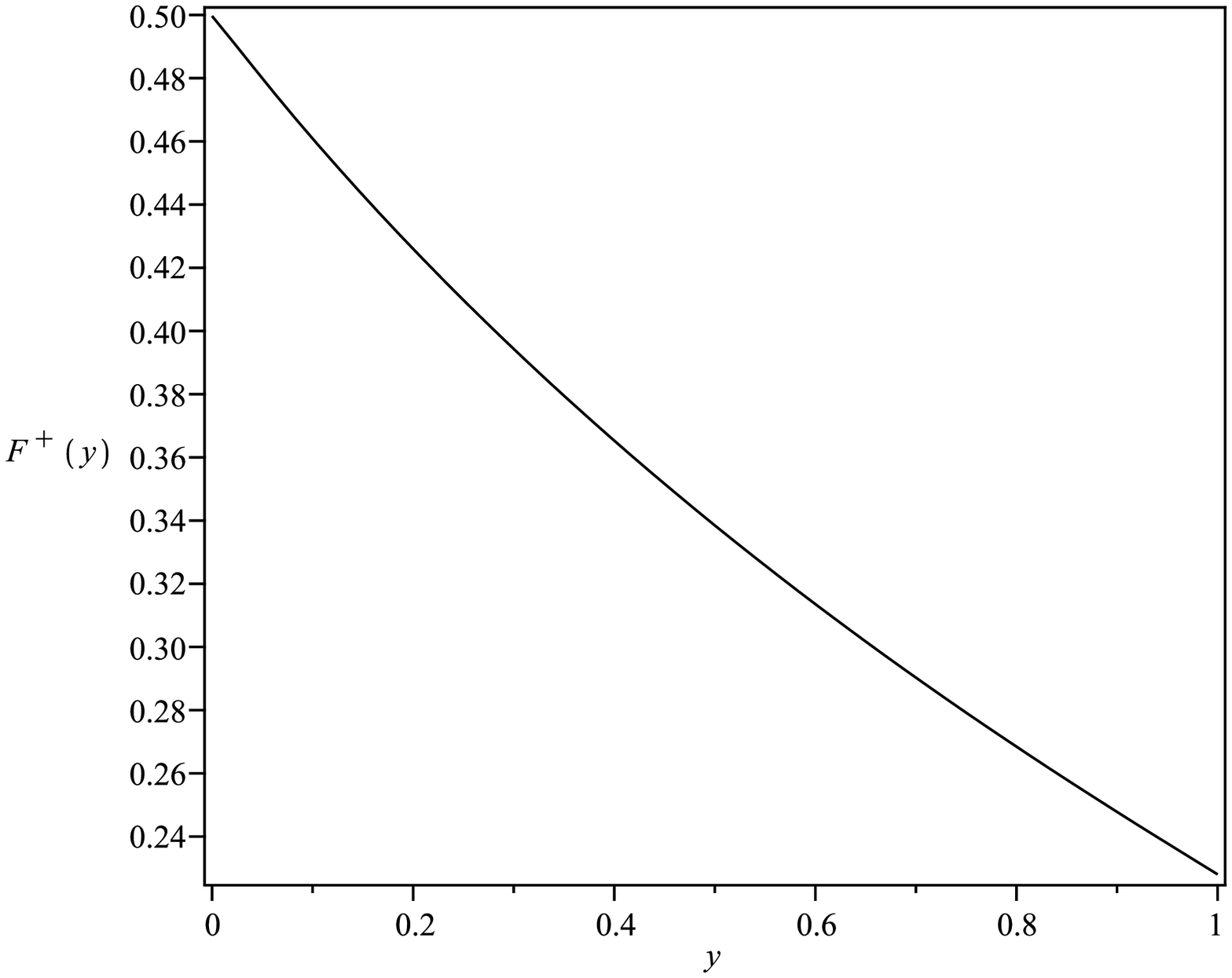}
\caption{Figure of $F^+(y)$ Example 2 by using MQ with $c=0.3$ and n=m=24.}
\label{Fig7}
\end{figure}

In this paper, a meshless method based on radial basis functions for solving the RTE in a slab medium with an isotropic scattering is
proposed. The recent approach solves the RTE numerically by using collocation nodes. RBFs are good ways to approximate multivariate
functions. They are proposed to provide an effective but simple way to improve the convergence of the solution by collocation method.
Now to show the accuracy of this method, we achieve $\|Res(y,x)\|^2$ by means of
\begin{eqnarray}
\|Res(y,x)\|^2=\int_{-1}^{1}\int_{0}^{1}Res^2(y,x)dy~dx~.
\end{eqnarray}
Table \ref{Tab4} shows values of $\|Res(y,x)\|^2$ for example $1$ by using Newton-Cotes integral approximate \cite{Abramowitz}.
It shows that by increasing $N$, value of $\|Res(y,x)\|^2$ decreases. Consequently, It provides to convergence of the method. \\
The resulting graphs of $Res(y,x)$ for examples $1,2$ are shown in Figures \ref{Fig2} and \ref{Fig3}.
Also the resulting graphs of $I(y,x)$ for example $1$ and example $2$ are shown in Figures \ref{Fig4} and \ref{Fig5}.
Figures \ref{Fig6} and \ref{Fig7} show the resulting graphs of $F^+(y)$.\\
Additionally, through the comparison, the exact values for radiative fluxes at
the lower boundary, with present method, and values of $\|Res(y,x)\|^2$, we have
showed that the RBFs approach has good reliability and efficiency.
\section*{Acknowledgements}
The research of first author (K. Parand) was supported by a grant from Shahid Beheshti University.

%*****************************************ACKNOWLEDGEMENTS**********************************************************************************
%\section*{Acknowledgements}
%The corresponding author would like to thank Shahid Beheshti University for the awarded grant.
%*********************************Bibliography***********************************************
\bibliographystyle{elsart}	 % (uses file "plain.bst")
\bibliography{myref}
%**************************************************Appendix***********************************************************************************
%**********************************************************Bibliography*************************************************
%\bibliographystyle{wileyj}	 % (uses file "plain.bst")
%\bibliography{ACM}
%%%%%%%%%%%%%%%%%%%%%%%%%%%%%%%%%%%%%%%%%%%%%%%%%%FIGURE%%%%%%%%%%%%%%%%%%%%%%%%%%%%%%%%%%%%%%%%%%%%%%%%%%%%%%%%%%%%
\clearpage
%&&&&&&&&&&&&&&&&&&&&&&&&&&&&&&&&&&&&&&&&&&&&&&&&&&&&&&&&&&&&&&&&&&&&&&&&&&&&&&&&&&&&&&&&&&&

%$$$$$$$$$$$$$$$$$$$$$$$$$$$$$$$$$$$$$$$$$$$$$$$$$$
%$$$$$$$$$$$$$$$$$$$$$$$$$$$$$$$$$$$$$$$$$$$$$$$$$$

%\clearpage
%........................
%.......................

%$$$$$$$$$$$$$$$$$$$$$$$$$$
%&&&&&&&&&&&&&&&&&&&&&&&&&&&&&&&&&&&&&&&&&&&&&&&&&&&&&&&&&&&&&&&&&&&&&&&&&&&&&&&&&&&&&&&&&&&

%%%%%%%%%%%%%%%%%%%%%%%%%%%%%%%%%%%%%%%%%%%%%%%%%%TABLE%%%%%%%%%%%%%%%%%%%%%%%%%%%%%%%%%%%%%%%%%%%%%%%%%%%%%%%%%%%%%
\clearpage
\begin{table}
\caption{Some well--known functions that generate RBFs $(r=\|x-x_{i}\|=r_{i}),~c>0$}
\begin{tabular}{l l}
\hline
%\multicolumn{2}{l}{}\\
%\cline{3-6}
Name of functions & Definition \\
%& Ref. \cite{Cheng and Pop} \\
\hline
Multiquadrics (MQ)  & $\sqrt{r^2+c^2}$\\
Inverse multiquadrics (IMQ) & $1/(\sqrt{r^2+c^2})$ \\
Gaussian (GA) & $exp(-cr^2)$ \\
Inverse quadrics & $1/(r^2+c^2)$ \\
\hline
\end{tabular}
\label{Tab.RBR.definition}
\end{table}
%&&&&&&&&&&&&&&&&&&&&&&&&&&&&&&&&&&&&&&&&&&table 1&&&&&&&&&&&&&&&&&&&&&&&&&&&&&&&&&&&&&&&&&&&&&&&&&&&&&&&&&&&
\begin{table}
\caption{Values of $F^+(1)$ for example $1$ with MQ, $m=n=20$ and $c=0.3$.}
\begin{tabular*}{\columnwidth}{@{\extracolsep{\fill}}*{9}{c}}
\hline
$c_{1}$&Algorithm&$t_{0}=0.1$&$t_{0}=0.5$&$t_{0}=1$&$t_{0}=3$\\
\hline
$0.7$ & GEA & Not reported & $0.753$    & $0.615$    & $0.369$\\
$$    & LSA & Not reported & $0.7498$   & $0.6112$   & $0.3547$\\
$$    & PLM & $0.93187$      & $0.75035$  & $0.61123$  & $0.35806$\\
$$    &Present method& $0.93071$     & $0.75049$ & $0.61211$ & $0.35834$\\
$$   &Exact & $0.931$        & $0.750$    & $0.611$    & $0.358$\\
\hline
$0.0$ & GEA & Not reported & $0.707$    & $0.555$    & $0.315$\\
$$    & LSA & Not reported & $0.7036$   & $0.5520$   & $0.2989$\\
$$    & PLM & $0.91710$      & $0.70434$  & $0.55340$  & $0.30131$\\
$$    &Present method& $0.91581$     & $0.70427$ & $0.55351$ & $0.30132$\\
$$   &Exact & $0.916$        & $0.704$    & $0.553$    & $0.301$\\
\hline
$-0.7$ & GEA &Not reported & $0.668$    & $0.507$    & $0.274$\\
$$    & LSA & Not reported & $0.6628$   & $0.5033$   & $0.2583$\\
$$    & PLM & $0.90242$      & $0.66327$  & $0.50483$  & $0.26007$\\
$$    &Present method& $0.901372$     & $0.663414$ & $0.504659$ & $0.260349$\\
$$   &Exact & $0.901$        & $0.663$    & $0.505$    & $0.260$\\
\hline
\end{tabular*}
\label{Tab2}
\end{table}
%\clearpage
%=====================================================================================
%....................................................
%...........................................................................................................................
\begin{table}
\caption{Value of $F^+(1)$ for example $2$ with MQ, $m=n=20$ and $c=0.3$.}
\begin{tabular*}{\columnwidth}{@{\extracolsep{\fill}}*{9}{c}}
\hline
Algorithm&$F^{+}(1)$&\vline&Algorithm&$F^{+}(1)$&$\|Res(y,x)\|^2$\\
\hline
GEA    & $0.458$ &\vline& PLM              & $0.4564$ & $-$\\
MTFM   & $0.471$ &\vline& Tau              & $0.4564$ & $-$\\
$P_{1}$  & $0.465$&\vline & Galerkin         & $0.4564$ & $-$\\
$P_{3}$  & $0.456$ &\vline&Present method$(n=10,m=10)$ & $0.457662$ & $9.4385e-04$\\
$F_{1}$  & $0.455$&\vline &Present method$(n=16,m=16)$ & $0.456551$ & $1.3043e-04$\\
$F_{3}$  & $0.456$ &\vline&Present method$(n=20,m=20)$ & $0.456254$ & $7.0662e-05$\\
$F_{9}$  & $0.456$ &\vline&Present method$(n=24,m=24)$ & $0.456096$ & $5.6480e-05$\\
\hline
\end{tabular*}
\label{Tab3}
\end{table}
%********************************************************************************8

%\clearpage
\begin{table}
\caption{Values of $\|Res(y,x)\|^2$ for example $1$ for case of $c_{1}=0.7$ with $c=0.3$.}
\begin{tabular*}{\columnwidth}{@{\extracolsep{\fill}}*{9}{c}}
\hline
$n=m$&Algorithm&$t_{0}=0.1$&$t_{0}=0.5$&$t_{0}=1$&$t_{0}=3$\\
\hline
$10$ & MQ & $5.2467e-03$ & $6.6703e-04$ & $1.5652e-04$ & $2.3255e-05$\\
$$   &IMQ & $4.0001e-02$ & $1.5039e-03$ & $5.3730e-04$ & $1.6458e-04$\\
$$   &IQ  & $2.7706e-01$ & $8.0117e-03$ & $2.3783e-03$ & $6.0969e-04$\\
\hline
$16$ &  MQ& $4.0359e-03$ & $3.0487e-04$ & $8.7236e-05$ & $6.9267e-06$\\
$$   &IMQ & $3.5955e-03$ & $2.7085e-04$ & $7.9158e-05$ & $7.3361e-06$\\
$$   & IQ & $4.5537e-03$ & $2.9087e-04$ & $8.6775e-05$ & $1.1738e-05$\\
\hline
$20$ & MQ & $3.5640e-03$ & $2.4073e-04$ & $7.1807e-05$ & $4.2414e-06$\\
$$   &IMQ & $3.0911e-03$ & $2.1409e-04$ & $6.4552e-05$ & $3.8916e-06$\\
$$   & IQ & $2.9274e-03$ & $2.0439e-04$ & $6.2031e-05$ & $4.0679e-06$\\
\hline
$24$ & MQ & $3.1107e-03$ & $2.0311e-04$ & $5.8004e-05$ & $2.4874e-06$\\
$$  & IMQ & $2.7001e-03$ & $1.8089e-04$ & $5.2018e-05$ & $2.1938e-06$\\
$$  &  IQ & $2.5609e-03$ & $1.7280e-04$ & $4.9889e-05$ & $2.1731e-06$\\
\hline
\end{tabular*}
\label{Tab4}
\end{table}
%\clearpage
%....................................................
%...........................................................................................................................
%...........................................................................................................................

%%%%%%%%%%%%%%%%%%%%%%%%%%%%%%%%%%%%%%%%%%%%%%%%%%FINISH%%%%%%%%%%%%%%%%%%%%%%%%%%%%%%%%%%%%%%%%%%%%%%
\end{document}